\documentclass[12pt]{amsart}
\usepackage[
  dvipdfmx,
  colorlinks=true,
  linkcolor=black,
  citecolor=black,
  urlcolor=blue
]{hyperref}
\usepackage{color}
\usepackage{amsmath,amssymb,amsthm}
\usepackage{amscd}
\usepackage{comment}
\usepackage{overpic}
\usepackage{tabularx}
\usepackage[dvipdfmx,left=25truemm,right=25truemm,top=30truemm,bottom=30truemm]{geometry}

\usepackage{bm}

\newtheorem{theorem}{Theorem}[section]
\newtheorem{lemma}[theorem]{Lemma}
\newtheorem{proposition}[theorem]{Proposition}

\theoremstyle{definition} 
\newtheorem{remark}[theorem]{Remark}

\newcommand{\R}{\mathbb{R}}

\renewcommand{\epsilon}{\varepsilon}

\renewcommand{\geq}{\geqslant}

\makeatletter
\@addtoreset{equation}{section}

\makeatother

\begin{document}

\title[On cylindrical directions of regular surfaces in $\R^4$]
{On cylindrical directions of regular surfaces in $\R^4$ derived from contact with spherical cylinders}

\author[M.~Hasegawa]
{Masaru Hasegawa}

\address[Masaru Hasegawa]{%
Department of Information Science, Center for Liberal Arts and Sciences, Iwate Medical University,
1-1-1 Idaidori, Yahaba-cho, Shiwa-gun, Iwate 028-3694, Japan.}
\email{mhase@iwate-med.ac.jp}


%


\allowdisplaybreaks[4]


\maketitle


\begin{abstract}
We investigate the contact of regular surfaces in $\mathbb{R}^4$ with spherical cylinders $S^2(r) \times \mathbb{R}$. 
We introduce cylindrical directions of these surfaces, defined by the axes of spherical cylinders having degenerate contact with the surfaces, and show that these directions capture the ridge condition of projected regular surfaces under orthogonal projection along these directions.
\end{abstract}


\section{Introduction}
\label{sec:introduction}

The geometric study of contact via singularity theory originated from Thom's insight into umbilic points and $D_4$ singularities (\cite{Thom1972}), which Porteous \cite{Porteous1971} later extended to characterize ridge curves via $A_3$-singularities of distance-squared functions. 
This correspondence was generalized by Montaldi \cite{Montaldi1986-1}, who proved that the geometric contact between arbitrary submanifolds can be characterized in terms of Mather's framework of $\mathcal{K}$-equivalence (\cite{Mather1968}).

The geometry of surfaces in $\R^4$ has been studied from the viewpoint of the above framework of contact and singularity theory.
Mochida, Romero-Fuster and Ruas \cite{MRFR1995} studied contact with hyperplanes, whereas Romero-Fuster and Sanabria-Codesal \cite{RFSC2013} studied contact with hyperspheres. 
Recently, Deolindo-Silva and Salarinoghabi \cite{DSS2026} studied contact with 2-planes. 

Takahashi \cite{Takahashi1970} showed a homogeneous surface in Euclidean $n$-space $\R^n$ is an orbit of a certain subgroup of the Euclidean motion group $O(n) \rtimes \R^n$ (see also \cite{Segre1938}). 
Geometrically, the homogeneous surfaces are surfaces which have constant principal curvatures at every point. 
In $\R^n$, the homogeneous surfaces are hyperplanes $\R^{n-1}$, hyperspheres $S^{n-1}(r)$ and hypercylinders $S^k(r) \times \R^{n-k-1}$.  
In particular, in $\R^3$, there are only three types homogeneous surfaces: planes, spheres or cylinders.  
Moreover, in $\R^4$, they are hyperplanes $\R^3$, hyperspheres $S^3(r)$, spherical cylinders $S^2(r) \times \R$ and toroidal cylinders $S^1(r) \times \R^2$. 

The contact of surfaces, including singular surfaces, in $\R^3$ with planes or spheres has been extensively studied (see, for example,  \cite{AOS2026, BGT1995, BGT1999, Francisco2023, FH2012-1, FH2012-2, FH2013, FH2025-1, FH2025-2, OSS2018, OST2018, Porteous1983, Teramoto2026}). 
In contrast, the contact of regular surfaces in $\R^3$ with cylinders had not been investigated in detail before the work of Fukui, the author and Nakagawa \cite{FHN2017}.   
They introduced the cylindrical direction of the surfaces via the contact with cylinders from the viewpoint of the singularity theory. 
A cylindrical direction $\bm{v}\in S^2$ at a point $\bm{p}$ of the surface is defined as the direction of the axis of a cylinder having degenerate contact with the surface, where the axis passes through the finite focal point of the regular apparent contour along $\bm{v}$, and the contour has a vertex at the point corresponding to $\bm{p}$ (see Theorems~\ref{thm:A_2} and \ref{thm:A_3} for details). 

Tari, Salarinoghabi and the author \cite{TSH2025} studied the deformations of apparent contours of regular surfaces in $\R^3$ using the fact that the cylindrical direction captures vertices of apparent contours of the surfaces. 
The author and Saji \cite{HS2025} attempted to extend the notion of cylindrical directions to cross caps and cuspidal edges in terms of contact with cylinders and the axes of cylinders having degenerate contact with the surfaces, which are related to vertices of their apparent contours. 
They showed that an analogue of cylindrical directions exists for cross caps, whereas the notion of cylindrical directions cannot be extended to cuspidal edges.

The aim of this paper is to investigate whether a higher-dimensional extension of the cylindrical direction of a regular surface in $\R^3$ exists through the contact of regular surfaces in $\R^4$ with spherical cylinders $S^2(r)\times\R$. 
Contact with spherical cylinders is related to orthogonal projections, as in the case of $\R^3$. 
Orthogonal projections of surfaces in $\R^4$ into 3-space have been studied from the viewpoint of singularity theory (see, for example, \cite{BROS2019, BN1998, DSilva2021, DSK2019, NBT2007, OST2015}). Recently, the author \cite{Hasegawa2026} studied the geometry of projected regular surfaces from the viewpoint of singularities of distance-squared functions on these surfaces.

Our main result is Theorem~\ref{thm:main}, which shows that the cylindrical direction of a regular surface in $\R^3$ can be extended to regular surfaces in $\R^4$ through contact with spherical cylinders. 
This higher-dimensional extension captures the ridge condition of a projected regular surface under the orthogonal projection along the axis direction of a spherical cylinder having degenerate contact with the original surface in $\R^4$.

The paper is organized as follows: In Section~\ref{sec:Preliminaries}, we recall the geometry of regular surfaces in $\R^4$, contact map-germs from the viewpoint of $\mathcal{K}$-equivalence, and the definition and geometric meaning of cylindrical directions of regular surfaces in $\R^3$.
In Section~\ref{sec:cylindrical_direction}, we give our main theorem and introduce a higher-dimensional extension of cylindrical directions. 
We also investigate a coefficient in the expression for the determinant of the Hessian matrix of the contact map-germ measuring contact between a regular surface in $\R^4$ and a spherical cylinder.

\section{Preliminaries}
\label{sec:Preliminaries}

\subsection{Differential geometry of surfaces in $\R^4$}
\label{sec:DGS}

While the primary focus of references such as \cite{BN1998, BT2002, GMRFR2000, Little1969, MRFR1995, NBT2007} lies elsewhere, the geometry of surfaces in $\R^4$ has been well documented therein. 
For a comprehensive discussion on the differential geometry of surfaces in $\R^4$ utilizing singularity theory, the reader is referred to \cite[Chapter 7]{IFRT2015}. 

The second-order extrinsic geometry of a regular surface in $\R^4$ is characterized by its vector-valued second fundamental form. 
Let $U$ be an open subset of $\R^2$ containing the origin. 
Consider a local parameterization of a regular surface $M \subset \R^4$ given by
$$
X\colon U \to \R^4
$$
such that $X(0,0)=\bm{p}\in M$. The set $\{X_u(0,0),X_v(0,0)\}$ then constitutes a tangent frame of $T_{\bm{p}}M$. 
Furthermore, let $\{\bm{\nu}_1,\bm{\nu}_2\}$ be an orthonormal frame of the normal space $N_{\bm{p}}M$. 
The \textit{first fundamental form at $\bm{p}$}, denoted by
$$
I_{\bm{p}}\colon T_{\bm{p}}M \times T_{\bm{p}}M \to \R,
$$
is defined via the relation
$$
I_{\bm{p}} = \langle dX, dX\rangle = E du^2+ 2F du dv+G dv^2,
$$
with the coefficients
$$
E=\langle X_u, X_u\rangle, \quad F=\langle X_u, X_v\rangle, \quad G=\langle X_v, X_v \rangle.
$$
The \textit{second fundamental form at $\bm{p}$} is a mapping
$$
II_{\bm{p}}\colon T_{\bm{p}}M \times T_{\bm{p}}M \to N_{\bm{p}}M,
$$
which is represented as
$$
II_{\bm p} = \langle\bm\nu_1,d^2 X\rangle\bm\nu_1 + \langle\bm\nu_2,d^2 X\rangle\bm\nu_2.
$$
Here, the component formulas are given by
$$
\langle\bm\nu_i,d^2 X\rangle = l_i\,du^2+2m_i\,du\,dv+n_i\,dv^2, \qquad i=1,2.
$$

Given a non-zero tangent direction $\bm{w} \in T_{\bm{p}}M\setminus\{\bm{0}\}$, the \textit{normal curvature vector} $\kappa_n$ at $\bm{p}$ is defined through
$$
\kappa_n(\bm{p},\bm{w}) = \dfrac{II_{\bm{p}}(\bm{w},\bm{w})}{I_{\bm{p}}(\bm{w},\bm{w})}.
$$
Since $\kappa_n(\bm{p})$ acts as a quadratic map, the image of the unit circle in $T_{\bm{p}}M$ forms an ellipse in the normal plane $N_{\bm{p}}M$, which may degenerate to a point or a radial line segment. 
This image is designated as the \textit{curvature ellipse at $\bm{p}$} and denoted by $\mathcal{E}_{\bm{p}}$. 
It follows that
\begin{equation*}
\mathcal{E}_{\bm{p}}=\{\kappa_n(\bm{p},\bm{w}): \bm{w}\in T_{\bm{p}}M, I_{\bm{p}}(\bm{w},\bm{w})=1\} \subset N_{\bm{p}}M.
\end{equation*}
The center of $\mathcal{E}_{\bm{p}}$ corresponds to the \textit{mean curvature vector $\bm{H}$ at $\bm{p}$}, which is explicitly expanded as
$$
\bm{H}=h_1\bm{\nu}_1+h_2\bm{\nu}_2=\dfrac{En_1-2Fm_1+Gl_1}{2(EG-F^2)}\bm{\nu}_1 + \dfrac{En_2-2Fm_2+Gl_2}{2(EG-F^2)}\bm{\nu}_2.
$$
The lengths of the two semi-axes of $\mathcal{E}_{\bm{p}}$ will be denoted by $\lambda$, $\mu\ge 0$.



In accordance with \cite{MRFR1995}, points $\bm{p}\in M$ are classified by the position of the origin $\bm{0}_{\bm{p}}$ of $N_{\bm{p}}M$ relative to a non-degenerate ellipse $\mathcal{E}_{\bm{p}}$.  
The point $\bm{p}$ is designated as \textit{elliptic}, \textit{parabolic}, or \textit{hyperbolic} if $\bm{0}_{\bm{p}}$ lies inside, on, or outside $\mathcal{E}_{\bm{p}}$, respectively. 
From this geometric definition, an elliptic, parabolic, or hyperbolic point $\bm{p}$ of $M$ is characterized by the existence of 0, 1, or 2 tangent lines to $\mathcal{E}_{\bm{p}}$ issuing from $\bm{0}_{\bm{p}}$, respectively. 
In the event that $\mathcal{E}_{\bm{p}}$ degenerates to a radial line segment, $\bm{p}$ is termed an \textit{inflection} point of $M$. 
Moreover, an inflection point $\bm{p}$ is said to be of real type, imaginary type, or flat type when $\bm{0}_{\bm{p}} \in \mathcal{E}_{\bm{p}}$, $\bm{0}_{\bm{p}}\not\in\mathcal{E}_{\bm{p}}$, or $\bm0_{\bm{p}}$ is an endpoint of $\mathcal{E}_{\bm{p}}$, respectively.

The curvature ellipse $\mathcal{E}_{\bm{p}}$ is viewed as the image of the unit circle in the tangent plane $T_{\bm{p}}M$ under a pair of quadratic forms $(Q_1,Q_2)$ associated with the second fundamental form $II_{\bm{p}}$.
By expressing a tangent vector as $\bm w=u\bm{e}_1+v\bm{e}_2$,
a pair of binary quadratic forms is determined by the second fundamental form as follows:
$$
Q_i(u,v)=l_i u^2+2m_i u v+n_i v^2, \quad i=1,2.
$$
The pencil of binary forms determined by $(Q_1,Q_2)$ was investigated in \cite{BN1998}.
Any binary quadratic form $A u^2+2B u v+C v^2$ can be identified with its coefficient vector $(A,B,C)\in\R^3$.
The cone of degenerate binary quadratic forms in $\R^3$ is defined by the equation $B^2-AC=0$.
A conic in $\R P^2$ is determined by this cone after projectivization.
Provided that $Q_1$ and $Q_2$ are linearly independent, a projective line in $\R P^2$ is determined by their pencil.
This line intersects the conic in zero, one or two points according as $\delta(\bm{p})$ is negative, zero, or positive, respectively, where
\begin{equation*}
\delta(\bm{p}) = (l_1 n_2-n_1 l_2)^2-4(l_1 m_2 - m_1 l_2)(m_1 n_2 - n_1 m_2)
\end{equation*}
A point $\bm{p} \in M$ is elliptic, parabolic, or hyperbolic according as $\delta(\bm{p})$ is negative, zero, or positive. 
Furthermore, $\bm p$ is an inflection point if the rank of $\mathcal A(\bm{p})$ equals 1, which implies that $Q_1$ and $Q_2$ are linearly dependent, where 
\begin{equation*}
\mathcal{A}(\bm{p}) = 
\begin{pmatrix}
l_1 & m_1 & n_1 \\
l_2 & m_2 & n_2
\end{pmatrix}
\end{equation*}
It is worth remarking that if $(Q_1,Q_2)$ is equivalent to $(x^2,0)$ or $(0,0)$ under the action of $GL(2,\R) \times GL(2,\R)$, then $\bm{p}$ is referred to as a \textit{degenerate inflection} point (see \cite{BT2002, NBT2007}). 
It should also be noted that the equivalence of $(Q_1,Q_2)$ to $(x^2,0)$ implies that $\mathcal{E}_{\bm{p}}$ is a radial line segment and $\bm{0}_{\bm{p}}$ is an endpoint of $\mathcal{E}_{\bm{p}}$, which corresponds to an inflection point of flat type. 
Conversely, if $(Q_1, Q_2)$ is equivalent to $(0,0)$, then $\mathcal{E}_{\bm{p}}$ collapses to a single point $\bm{0}_{\bm{p}}$, which means $\bm{p}$ is not an inflection point in the sense of \cite{MRFR1995}. 


Let $\bm{w}(\theta)\in S(T_{\bm{p}}M)$, $\theta\in[0,2\pi)$, be a regular parameterization of the unit circle in $T_{\bm{p}}M$.
An \textit{asymptotic direction} of $M$ at $\bm{p}$ is defined as a direction $\bm{w}(\theta_0) \in T_{\bm{p}}M$ along which 
$$
\kappa_n(\bm{p},\bm{w}(\theta_0))\quad \text{and} \quad \dfrac{\partial}{\partial \theta}\kappa_n(\bm{p},\bm{w}(\theta_0))
$$
are parallel, where $\bm{w}(\theta)$ and $-\bm{w}(\theta)$ are identified as the same direction (see \cite{GMRFR2000, MRFR1995}).  
Equivalently, the status of $\bm{w}(\theta_0)$ as an asymptotic direction is achieved if the tangent line to $\mathcal{E}_{\bm{p}}$ at $\kappa_n(\bm{p},\bm{w}(\theta_0))$ passes through $\bm{0}_{\bm{p}}$ (see \cite{NBT2007}).
Therefore, under the assumption that $\bm{p}$ is neither an inflection point nor a degenerate inflection point, the number of asymptotic directions at $\bm{p}$ is zero, one or two, depending on whether $\bm{p}$ is an elliptic, parabolic, or hyperbolic point of $M$, respectively. 
Every direction in $T_{\bm{p}}M$ becomes asymptotic if $\bm{p}$ is an inflection point (see \cite{MRFR1995}). 
In \cite{BT2002}, the binary differential equation 
\begin{equation}
\label{eq:BDE_of_asymptotic}
(l_1 m_2-m_1l_2)du^2+(l_1 n_2-n_1l_2)du dv+(m_1n_2-n_1m_2)dv^2=0
\end{equation}
was introduced, and the asymptotic directions of $M$ at $\bm{p}$ were shown to be solutions of \eqref{eq:BDE_of_asymptotic}. 
The discriminant of \eqref{eq:BDE_of_asymptotic} is identical to $\delta(\bm{p})$. 

For a unit normal vector $\bm\nu=\cos\phi \bm{\nu}_1+\sin\phi\bm{\nu}_2 \in S(N_{\bm{p}}M)$, $\phi\in(-\pi/2,\pi/2]$, a symmetric bilinear form 
$$
II_{\bm{p},\bm\nu} \colon T_{\bm{p}}M \times T_{\bm{p}}M \to \R \quad\text{is defined by}\quad II_{\bm{p},\bm\nu}(\bm{w}_1,\bm{w}_2)=\langle II_{\bm{p}}(\bm{w}_1,\bm{w}_2), \bm\nu \rangle,\quad (\bm{w}_i \in T_{\bm{p}}M)
$$
Hence, in terms of the tangent frame $\{X_u,X_v\}$ of $T_{\bm{p}}M$, we obtain 
\begin{equation*}
\begin{split}
\left[II_{\bm{p},\bm\nu}\right] &= 
\begin{pmatrix}
II_{\bm{p},\bm\nu}(X_u,X_u) & II_{\bm{p},\bm\nu}(X_u,X_v) \\
II_{\bm{p},\bm\nu}(X_u,X_v) & II_{\bm{p},\bm\nu}(X_v,X_v)
\end{pmatrix} \\
&=
\begin{pmatrix}
l_1\cos\phi+l_2\sin\phi & m_1\cos\phi+m_2\sin\phi \\
m_1\cos\phi+m_2\sin\phi & n_1\cos\phi+n_2\sin\phi
\end{pmatrix}.
\end{split}
\end{equation*}
A direct computation yields 
\begin{equation*}
\text{det}[II_{\bm{p},\bm\nu}] = (l_1 n_1-m_1^2)\cos^2\phi+(l_1n_2-2m_1m_2+n_1l_2)\cos\phi\sin\phi+(l_2n_2-m_2^2)\sin^2\phi.
\end{equation*}
Moreover, by treating $\det[II_{\bm p,\bm\nu}]$ as a homogeneous quadratic form in $(\cos\phi,\sin\phi)$, its discriminant is seen to be
\begin{equation*}
\text{disc}(\text{det}[II_{\bm{p},\bm\nu}])= \delta(\bm{p}).
\end{equation*}


The condition $\text{det}[II_{\bm{p},\bm\nu}]=0$ is equivalent to the existence of a non-zero vector $\bm{w}_0 \in T_{\bm{p}}M$ satisfying $II_{\bm{p},\bm\nu}(\bm{w}_0,\bm{w})=0$ for every $\bm{w} \in T_{\bm{p}}M$. 
Therefore, regarding $\bm{w}_0$ as the kernel direction of $[II_{\bm{p},\bm\nu}]$, the relation $\langle II_{\bm{p}}(\bm{w}_0,\bm{w}_0),\bm\nu \rangle=0$ is obtained, which implies $\langle II_{\bm{p}}(\bm{w}_0,\bm{w}_0^{\perp}),\bm\nu \rangle=0$. 
Fix a regular parameterization $\bm{w}(\theta)\in S(T_{\bm{p}}M)$, $\theta\in[0,2\pi)$, of the unit circle in $T_{\bm{p}}M$.
We then obtain
\begin{equation} 
\label{eq:kappa_and_kappa_theta} \kappa_n(\bm{p}, \bm{w}(\theta))=II_{\bm{p}}(\bm{w}(\theta),\bm{w}(\theta))\quad \text{and} \quad \dfrac{\partial}{\partial \theta} \kappa_n(\bm{p},\bm{w}(\theta))=2a(\theta)II_{\bm{p}}(\bm{w}(\theta),\bm{w}(\theta)^{\perp}), \end{equation}
where 
$$ 
a(\theta)\ne0\quad \text{with} \quad \frac{\partial}{\partial \theta}\bm{w}(\theta) = a(\theta)\bm{w}(\theta)^{\perp}. 
$$
Evaluated at the kernel direction $\bm{w}_0=\bm{w}(\theta_0)$, equations \eqref{eq:kappa_and_kappa_theta} yield
\begin{equation} 
\label{eq:II_p_nu} 
\begin{split} II_{\bm{p},\bm\nu}(\bm{w}_0,\bm{w}_0) &= \langle\kappa_n(\bm{p},\bm{w}(\theta_0)),\bm\nu \rangle=0, \\ II_{\bm{p},\bm\nu}(\bm{w}_0,\bm{w}_0^{\perp}) & = \langle II_{\bm{p}}(\bm{w}_0,\bm{w}_0^{\perp}),\bm\nu\rangle = \dfrac1{2a(\theta_0)}\left\langle \dfrac{\partial}{\partial \theta}\kappa_n(\bm{p},\bm{w}(\theta_0)),\bm\nu \right\rangle =0. \end{split} 
\end{equation}
Equation \eqref{eq:II_p_nu} immediately implies that both
$$ 
\kappa_n(\bm{p},\bm{w}(\theta_0)) \quad\text{and}\quad \frac{\partial}{\partial\theta}\kappa_n(\bm{p},\bm{w}(\theta_0)) 
$$
are contained in the one-dimensional subspace $\bm\nu^\perp$ of $N_{\bm{p}}M$.
Their parallelism ensures that the kernel direction $\bm{w}_0$ of $[II_{\bm{p},\bm\nu}]$ serves as an asymptotic direction at $\bm{p}$.
It is worth noting that at an inflection point $\bm{p}$, or at a point $\bm{p}$ where the curvature ellipse degenerates to a point distinct from $\bm{0}_{\bm{p}}$, one can find a unit normal vector $\bm\nu_0$ satisfying $[II_{\bm{p},\bm\nu_0}]=\bm{0}$.
On the other hand, if the curvature ellipse degenerates to $\bm{0}_{\bm{p}}$ at $\bm{p}$, then $[II_{\bm{p},\bm\nu}]=\bm{0}$ holds for any unit normal vector $\bm\nu$.
As a result, every tangent direction at $\bm{p}$ is asymptotic under these conditions.

Conversely, assume $\bm{w}_0=\bm{w}(\theta_0)$ to be an asymptotic direction at $\bm{p}$.
The tangent line $\ell$ to $\mathcal{E}_{\bm{p}}$ at $\kappa_n(\bm{p},\bm{w}(\theta_0))$ thus passes through $\bm{0}_{\bm{p}}$.
Let $\bm\nu\in N_{\bm{p}}M$ be a unit normal vector to $\ell$.
The fact that $\ell$ passes through $\bm{0}_{\bm{p}}$ allows us to identify $\ell$ with the one-dimensional subspace $\bm\nu^\perp$ of $N_{\bm{p}}M$.
From \eqref{eq:kappa_and_kappa_theta}, it follows that
$$
 II_{\bm{p}}(\bm{w}_0,\bm{w}_0) =\kappa_n(\bm{p},\bm{w}(\theta_0))\in\bm\nu^\perp \quad \text{and}\quad II_{\bm{p}}(\bm{w}_0,\bm{w}_0^\perp) =\frac{1}{2a(\theta_0)} \frac{\partial}{\partial\theta} \kappa_n(\bm{p},\bm{w}(\theta_0)) \in\bm\nu^\perp. 
$$
Consequently,
$$ 
\langle II_{\bm{p}}(\bm{w}_0,\bm{w}_0),\bm\nu\rangle = II_{\bm{p},\bm\nu}(\bm{w}_0,\bm{w}_0) = 0\quad \text{and} \quad \langle II_{\bm{p}}(\bm{w}_0,\bm{w}_0^{\perp}),\bm\nu \rangle = II_{\bm{p},\bm\nu}(\bm{w}_0,\bm{w}_0^\perp) =0. 
$$
Given that $\{\bm{w}_0,\bm{w}_0^\perp\}$ forms an orthonormal basis of $T_{\bm{p}}M$, we conclude that $II_{\bm{p},\bm\nu}(\bm{w}_0,\bm{w})=0$ for every $\bm{w}\in T_{\bm{p}}M$.
Hence, $\bm{w}_0$ is the kernel direction of $[II_{\bm{p},\bm\nu}]$.
This normal direction $\bm\nu$, regarded up to sign, is called a \textit{binormal direction} of $M$ at $\bm{p}$.
While exactly two binormal directions exist at a hyperbolic point, they coincide at a parabolic point (cf. \cite{MRFR1995,NBT2007}).
Equivalently, a normal direction $\bm\nu$ is binormal if and only if the height function in the direction $\bm\nu$ exhibits a degenerate critical point at $\bm{p}$.
This characterization is also applicable at inflection points (see \cite{MRFR1995}) and at points where the curvature ellipse degenerates to a point, even if a unique asymptotic direction cannot be determined.



In \cite{Hasegawa2026}, the author introduced a curvature-ellipse-based Monge normal form of a regular surface in $\R^4$, which is normalized by only Euclidean motions of the ambient space and changes of coordinates in the source of the parameterization. 
Therefore, the normal form preserves the Euclidean geometry of the surface, specifically, the center and semi-axis lengths of the curvature ellipse up to their ordering. 


Consequently, for the study of the geometry of the surface, including analysis of contact between the surface with model submanifolds, we use this normal form instead of a normal form derived from a general $GL(2,\R)\times GL(2,\R)$ action applied to pairs of binary forms. 
Even though such a normal form is useful for the study of affine or flat geometry of the surface, it does not generically preserve the Euclidean geometry. 
 
In that sense, to investigate the contact between surfaces in $\R^4$ and spherical cylinders, we therefore adopt the normal form given by Proposition~\ref{prop:Euclidean_Monge_normal_form}, including the cases for which the curvature ellipse degenerates.

%

\begin{proposition}[{\cite[Theorem 2.1]{Hasegawa2026}}]
\label{prop:Euclidean_Monge_normal_form}
Let $M$ be a regular surface in $\R^4$, and let $\bm{p}\in M$.
After a Euclidean motion of $\R^4$ and a local reparameterization of $M$, there exists a local Monge normal form at $\bm{p}$ such that
\begin{equation*}
X(u,v)=\left(u,v,f(u,v),g(u,v)\right),
\end{equation*}
where
\begin{align*}
f(u,v) &= \frac{1}{2}\left\{h_1(u^2+v^2)+\lambda(u^2-v^2) \right\}+\sum_{n\ge3}f_n(u,v),\quad f_n(u,v)=\sum_{i+j=n}\dfrac1{i!j!}a_{ij}u^iv^j, \\
g(u,v) &= \frac{1}{2} \left\{h_2(u^2+v^2)+2\mu uv\right\}+\sum_{n\ge3}g_n(u,v),\quad g_n(u,v)=\sum_{i+j=n}\dfrac1{i!j!}b_{ij}u^iv^j.
\end{align*}
\end{proposition}

\subsection{Contact}
\label{sec:contact}

We begin with the notion of contact type introduced by Montaldi \cite{Montaldi1986-1} in the framework of singularity theory.

Let $(N_1,Y_1)$ and $(N_2,Y_2)$ be two pairs of submanifold-germs at the origin of $\R^n$.
The pairs are said to {\it have the same contact type} if there exists a diffeomorphism-germ $\alpha \colon (\R^n, \bm{0}) \to (\R^n, \bm{0})$ such that $\alpha(N_1)=N_2$ and $\alpha(Y_1)=Y_2$.

The following theorem relates contact type to $\mathcal{K}$-equivalence. 
Here, we say that two map-germs $f_1, f_2 \colon (\R^n, \bm{0}) \to (\R^p,\bm{0})$ are {\it $\mathcal{K}$-equivalent} if there exist a diffeomorphism-germ $\alpha \colon (\R^n, \bm{0}) \to (\R^n, \bm{0})$ and a smooth map $A \colon (\R^n, \bm{0}) \to GL(p,\R)$ such that $A f_1=f_2\circ \alpha.$

\begin{theorem}[\cite{Montaldi1986-1}]
\label{thm:contact}
Let $\iota_i \colon (N_i,x_i) \to (\R^n, \bm{0})$ $(i=1,2)$ be immersion-germs and let $\beta_i\colon (\R^n, \bm{0}) \to (\R^p, \bm{0})$ be submersion-germs with $Y_i=\beta_i^{-1}(\bm{0})$.
Then the pairs $(N_1,Y_1)$ and $(N_2,Y_2)$ have the same contact type if and only if $\beta_1 \circ \iota_1$ and $\beta_2 \circ \iota_2$ are $\mathcal{K}$-equivalent.
\end{theorem}

Theorem~\ref{thm:contact} implies that if $\iota \colon (N,x) \to (\R^n,\bm{0})$ is an immersion-germ and $\beta \colon (\R^n, \bm{0}) \to (\R^p, \bm{0})$ is a submersion-germ with $Y=\beta^{-1}(\bm{0})$, then the contact type of $N$ and $Y$ is determined by the $\mathcal{K}$-equivalence class of the map-germ $\beta \circ \iota$. 
We call the map-germ the {\it contact map-germ} of $N$ and $Y$.

In this paper, we consider the following $\mathcal{K}$-singularities (see \cite{Arnold1986}):
$$
A_k^{\pm}: x^2 \pm  y^{k+1} \quad (k \geq 1), \quad D_k^{\pm}: x^2y \pm y^{k-1} \quad (k \geq 4). 
$$


\subsection{Cylindrical directions of regular surfaces in $\R^3$}
\label{sec:cylindrical_direction}

In \cite{FHN2017}, the authors introduced the notion of the cylindrical directions of regular surfaces $M_3$ in $\R^3$ via their contact with cylinders.

Let $C^{(3)}_{\bm{v},\bm{a},r}$ denote a cylinder whose axis is parallel to $\bm{v} \in S^2$, passes through $\bm{a} \in \bm{v}^{\perp}$ but does not pass through $\bm{p} \in M_3$ and whose cross-sectional circle has of radius $|r| \ne0$. 

Assume that $\bm{p}$ is a common point of $M_3$ and $C^{(3)}_{\bm{v},\bm{a},r}$. 
Then $M_3$ and $C^{(3)}_{\bm{v},\bm{a},r}$ are not tangent at $\bm{p}$ when they do not have a common normal line there. 
So they have $A_{\geq1}$-contact at $\bm{p}$ if and only if $M_3$ and $C^{(3)}_{\bm{v},\bm{a},r}$ have a common normal line at $\bm{p}$. 
This condition is equivalent to that $\bm{a}$ is on the projected line of the common normal line there under the orthogonal projection along $\bm{v}$ onto $\bm{v}^{\perp}$. 
In this case, $\bm{v} \in T_p M_3$, and hence the following theorems hold:

\begin{theorem}[{\cite[Theorem 4.1]{FHN2017}}]
\label{thm:A_2}
The contact type of $M_3$ with $C^{(3)}_{\bm{v},\bm{a},r}$ at $\bm{p}$ is $A_{\geq 2}$-contact if and only if one of the following conditions holds: 
\begin{enumerate}
\item
$\bm{p}$ is not a parabolic point, $\bm{v}$ is not an asymptotic direction at $\bm{p}$ and $\bm{a}$ is the focal point of the apparent contour of $M_3$ along $\bm{v}$ at a point corresponding to $\bm{p}$. 

\item
$\bm{p}$ is a parabolic point but not an umbilic point, $\bm{v}$ is an asymptotic direction at $\bm{p}$ and $\bm{a}$ is not a point corresponding to the focal point at $\bm{p}$. 

\item
$\bm{p}$ is a flat umbilic point. 
\end{enumerate}
\end{theorem}

\begin{theorem}[{\cite[Proposition 4.12]{FHN2017}}]
\label{thm:A_3}
If Theorem~\ref{thm:A_2}~(1) holds, then the contact type of $M_3$ with $C^{(3)}_{\bm{v},\bm{a},r}$ at $\bm{p}$ is $A_{\geq k}$-contact $(k \geq 3)$ if and only if the apparent contour has a $(k-2)$-order vertex at a point corresponding to $\bm{p}$.
\end{theorem}

We say that the direction $\bm{v}$ of the axis of $C^{(3)}_{\bm{v},\bm{a},r}$ is a \textit{cylindrical direction of $M_3$ at $\bm{p}$} if the contact type of $M_3$ with $C^{(3)}_{\bm{v},\bm{a},r}$ at $\bm{p}$ is $A_{\geq 3}$-contact under the assumption of Theorem~\ref{thm:A_2} (1).

We remark that the notion of a cylindrical direction is not determined solely by the $A_{\geq 3}$-contact condition, since the cases in Theorem~\ref{thm:A_2} (2) and (3) are excluded from the definition.
More precisely, this definition requires that the apparent contour of $M$ along $\bm{v}$ be a regular curve and have not inflection point at the point corresponding to $\bm{p}$, that this point be a vertex of the apparent contour, and that the axis of $C^{(3)}_{\bm{v},\bm{a},r}$ pass through its finite focal point.


\section{Cylindrical directions of regular surfaces in $\R^4$}
\label{sec:cylindrical_direction}

In \cite{FHN2017}, the authors analyzed the contact of a regular surface in $\R^3$ with cylinders by using the contact map-germ defined by the difference of third components of a Monge form of the surface and the cylinder.  
However, in this paper, directly based on the concept of Montaldi's contact map-germ, we investigate the contact map-germ defined by the composite of a parameterization of a regular surface in $\R^4$ with a defining function whose zero set defines the spherical cylinder $S^2(r) \times \R$ (cf. \cite{HS2025}). 

Assume that $X \colon (\R^2,\bm{q}) \to (\R^4,\bm{p})$ is a smooth map-germ parameterizing a regular surface $M$ in $\R^4$.
Let denote $C_{\bm{v},\bm{a},r}$ a spherical cylinder whose radius is $|r|$, where $r\ne0$, and axis is parallel to $\bm{v}\in S^3$, passes through $\bm{a}\in \bm{v}^{\perp}$ but does not pass through $\bm{p}$. 
Hence, $C_{\bm{v},\bm{a},r}$ is given by
$$
C_{\bm{v},\bm{a},r} =\left\{\bm{x}\in\R^4 \colon \lVert \bm{x} - \langle \bm{x},\bm{v} \rangle \bm{v} - \bm{a} \rVert^2 = r^2 \right\}.
$$
Let $\bm{p}$ lie in the intersection between $M$ and $C_{\bm{v},\bm{a},r}$. 
We define the spherical-cylinder-contact map-germ associated to $C_{\bm{v},\bm{a},r}$ by
$$
F_C(u,v) = \lVert X(u,v) - \langle X(u,v),\bm{v} \rangle \bm{v} - \bm{a} \rVert^2. 
$$

Let 
$$
\pi_{\bm{v}}=P_{\bm{v}} \circ X \colon (\R^2,\bm{q}) \to \bm{v}^{\perp}, \quad \pi_{\bm{v}}(u,v) = X(u,v)-\langle X(u,v),\bm{v}\rangle \bm{v}
$$
be the orthogonal projection of $M$ onto $\bm{v}^{\perp}$, where $P_{\bm{v}}$ denotes the orthogonal projection onto $\bm{v}^{\perp}$:
$$
P_{\bm{v}}\colon \R^4 \to \bm{v}^{\perp},\quad P_{\bm{v}}(\bm{x})=\bm{x}-\langle \bm{x},\bm{v} \rangle \bm{v}.
$$

We denote $P_{\bm{v}}(M)$ by $\widetilde{M}$ the projected surface. 
The projection $\pi_{\bm{v}}$ is an immersion at $\bm{q}$ if and only if $\bm{v}\not\in T_{\bm{p}}M$. 
In this case, $\widetilde{M}$ is a regular surface at $P_{\bm{v}}(\bm{p})$. 
On the other hand, the projection $\pi_{\bm{v}}$ is singular at $\bm{q}$ if and only if $\bm{v}\in T_{\bm{p}}M$. 
In this case, $\widetilde{M}$ is a singular surface at $P_{\bm{v}}(\bm{p})$ (see, for example, \cite{NBT2007}). 

The family of projections is given by
$$
\Pi\colon M \times S^3 \to TS^3,\quad \Pi(\bm{p},\bm{v})=(\bm{v}, \bm{p}-\langle \bm{p},\bm{v}\rangle \bm{v}).
$$
For a fixed direction $\bm{v}$, the projection can be considered locally at a point $\bm{p}\in M$ as a map-germ 
$$
\pi_{\bm{v}}\colon (\R^2,\bm{0}) \to (\R^3,\bm{0}).
$$
The generic singularities of $\pi_{\bm{v}}$ are classified under the $\mathcal{A}$-equivalence. 
Here, two map-germs $f_1, f_2\colon(\R^2,\bm{0})\to(\R^3,\bm{0})$ are said to be \textit{$\mathcal{A}$-equivalent} if $f_2=\Psi \circ f_1 \circ \Phi^{-1}$ for some diffeomorphisms-germs $\Phi$ and $\Psi$ of the source and target, respectively. 
%
%

The study of singularities of orthogonal projections of surfaces in $\R^4$ into $3$-spaces goes back to Mond \cite{Mond1982, Mond1985}. 
He adopted a modification of the transversality theorem of David \cite{David1983}. 
Following these studies, subsequent studies have further developed this topic (see, for example, \cite{BROS2019, BN1998, DSilva2021, DSK2019, NBT2007, OST2015}). 
However, these studies have mainly focused on singularities of projections and the geometry of projected singular surfaces, although some aspects of the geometry of projected regular surfaces were also considered in \cite{NBT2007}. 
More recently, the author \cite{Hasegawa2026} studied the geometry of projected regular surfaces from the viewpoint of singularities of distance-squared functions on these surfaces.


We first determine the condition for $F_C$ to have a singular point.

\begin{proposition}
\label{prop:A_1}
The map-germ $F_C$ has an $A_{\geq1}$-singularity at $\bm{q}$ if and only if there exists a unit normal vector $\bm\nu \in N_{\bm{p}}M \cap \bm{v}^{\perp}$ such that $\bm{a}=P_{\bm{v}}(\bm{p})+r\bm\nu$. 
\end{proposition}
\begin{proof}
Using $\lVert\bm{v}\rVert=1$ and $\langle\bm{a},\bm{v}\rangle=0$, we obtain
$$
(F_C)_u(\bm{q}) = 2\langle X_u(\bm{q}), \pi_{\bm{v}}(\bm{q}) - \bm{a} \rangle, \quad
(F_C)_v(\bm{q}) = 2\langle X_v(\bm{q}), \pi_{\bm{v}}(\bm{q})-\bm{a}\rangle.
$$
Since $\{X_u(\bm{q}),X_v(\bm{q})\}$ is a basis of $T_{\bm{p}}M$, $(F_C)_u(\bm{q})=(F_C)_v(\bm{q})=0$ if and only if 
$$
\bm{a}=\pi_{\bm{v}}(\bm{q})+r\bm\nu = P_{\bm{v}}(\bm{p})+r\bm\nu
$$
for a suitable choice of $r\ne0$ and a unit normal vector $\bm\nu \in N_{\bm{p}}M \cap \bm{v}^{\perp}$. 
\end{proof}

\begin{remark}
The condition that $F_C$ has an $A_{\geq1}$-singularity at $\bm{q}$ does not imply that $\bm{v}\in T_{\bm{p}}M$. 
Indeed, if $\bm\nu' \in N_{\bm{p}}M$ is a unit normal vector orthogonal to $\bm\nu$, then
$$
\bm\nu^\perp = T_{\bm{p}}M \oplus \text{span}\{\bm\nu'\}.
$$
Thus, the axis direction $\bm{v}$ may have a non-zero component in the
$\bm\nu'$-direction. 
In particular, unlike the case of surfaces in $\R^3$, the distinction between tangential and non-tangential projection does not yet arise from the $A_{\geq 1}$-contact condition.
\end{remark}

To characterize degenerate singularities of the contact map-germ $F_C$ in terms of the geometry of projected regular surfaces $\widetilde{M}$, we recall the following result on the distance-squared function-germ  $d_{\bm{a}}$ on $\widetilde{M}$ from a point $\bm{a} \in \bm{v}^{\perp} \subset \R^4$ defined by
$$
d_{\bm{a}}\colon (\R^2,\bm{q}) \to \R, \quad d_{\bm{a}}(u,v)=\lVert \pi_{\bm{v}}(u,v)-\bm{a}\rVert^2.
$$

\begin{lemma}[{\cite[Proposition 4.3--4.5]{Hasegawa2026}}]
\label{lem:main}
Assume that $\widetilde{M}$ is a regular surface at $P_{\bm{v}}(\bm{p})$. 
\begin{enumerate}
\item
$d_{\bm{a}}$ has an $A_{\ge2}$-singularity at $\bm{q}$ if and only if $P_{\bm{v}}(\bm{p})$ is a non-umbilic point of $\widetilde{M}$ and, for some $i=1,2$,  
\begin{equation*}
\begin{array}{cccccc}
 & \bm{v}^\perp & \subset & \R^4 & \\ [-1mm]
 & \rotatebox[origin=c]{90}{$\in$}
 & & \rotatebox[origin=c]{90}{$\in$} & \\ [-1mm]
\bm{a}= & P_{\bm{v}}(\bm{p})+\dfrac1{\widetilde{\kappa}_i(\bm{q})}\widetilde{\bm{n}}(\bm{q})
& = &P_{\bm{v}}(\bm{p})+\dfrac1{\widetilde{\kappa}_i(\bm{q})} \bm\nu & (\widetilde{\kappa}_i(\bm{q})\ne0),
\end{array}
\end{equation*}
where $\widetilde{\kappa}_i$ is the principal curvature of $\widetilde{M}$ and $\widetilde{\bm{n}}$ is the unit normal vector field to $\widetilde{M}$.

\item
$d_{\bm{a}}$ has a $D_{\ge4}$-singularity at $\bm{q}$ if and only if 
$$
\bm{a}=P_{\bm{v}}(\bm{p})+\dfrac1{\widetilde{\kappa}(\bm{q})}\widetilde{\bm{n}}(\bm{q})= P_{\bm{v}}(\bm{p})+\dfrac1{\widetilde{\kappa}(\bm{q})} \bm\nu
$$ 
and $P_{\bm{v}}(\bm{p})$ is a non-flat umbilic point of $\widetilde{M}$, where $\tilde{\kappa}$ is the common principal curvature.

\item
$d_{\bm{a}}$ has an $A_{\ge3}$-singularity at $\bm{q}$ if and only if $P_{\bm{v}}(\bm{p})$ is a non-umbilic point of $\widetilde{M}$, for some $i=1,2$, 
$$
\bm{a}=P_{\bm{v}}(\bm{p})+\dfrac1{\widetilde{\kappa}_i(\bm{q})}\widetilde{\bm{n}}(\bm{q}) = P_{\bm{v}}(\bm{p})+\dfrac1{\widetilde{\kappa}_i(\bm{q})} \bm\nu \quad (\widetilde{\kappa}_i(\bm{q})\ne0)
$$
and $P_{\bm{v}}(\bm{p})$ is a ridge point relative to $\widetilde{\bm{v}}_i$, where $\widetilde{\bm{v}}_i$ is the principal direction with respect to $\widetilde{\kappa}_i$.
\end{enumerate}
\end{lemma}
The reader is referred to \cite{Hasegawa2026} for the condition for $d_{\bm{a}}$ having each singularity in terms of the geometry of $M$ and to \cite{Porteous1971, Porteous1983} (or \cite{FH2012-1} for another example) for the characterization of singularities of distance-squared functions on regular surfaces in $\R^3$.

The following theorem gives a characterization of degenerate singularities of $F_C$ by the geometry of projected regular surfaces $\widetilde{M}$.
 
\begin{theorem}
\label{thm:main}
Assume that $\widetilde{M}$ is a regular surface at $P_{\bm{v}}(\bm{p})$. 
\begin{enumerate}
\item
$F_C$ has an $A_{\geq 2}$-singularity at $\bm{q}$ if and only if $P_{\bm{v}}(\bm{p})$ is a non-umbilic point and for some $i=1,2$ 
$$
\bm{a} = P_{\bm{v}}(\bm{p})+\dfrac1{\widetilde{\kappa}_i(\bm{q})} \widetilde{\bm{n}}(\bm{q}) = P_{\bm{v}}(\bm{p})+\dfrac1{\widetilde{\kappa}_i(\bm{q})} \bm\nu \quad (\widetilde{\kappa}_i(\bm{q})\ne0).
$$ 

\item
$F_C$ has a $D_{\geq 4}$-singularity at $\bm{q}$ if and only if $P_{\bm{v}}(\bm{p})$ is a non-flat umbilic point and 
$$
\bm{a} = P_{\bm{v}}(\bm{p})+\dfrac1{\widetilde{\kappa}(\bm{q})} \widetilde{\bm{n}}(\bm{q})=P_{\bm{v}}(\bm{p})+\dfrac1{\widetilde{\kappa}(\bm{q})} \bm\nu.
$$

\item
$F_C$ has an $A_{\geq 3}$-singularity at $\bm{q}$ if and only if $P_{\bm{v}}(\bm{p})$ is a non-umbilic point, for some $i=1,2$ 
$$
\bm{a} = P_{\bm{v}}(\bm{p})+\dfrac1{\widetilde{\kappa}_i(\bm{q})} \widetilde{\bm{n}}(\bm{q}) = P_{\bm{v}}(\bm{p})+\dfrac1{\widetilde{\kappa}_i(\bm{q})} \bm\nu \quad (\widetilde{\kappa}_i(\bm{q})\ne0)
$$
and $P_{\bm{v}}(\bm{p})$ is a ridge point relative to $\widetilde{\bm{v}}_i$.
\end{enumerate} 
\end{theorem}

\begin{proof}
By applying a suitable Euclidean motion, we may assume without loss of generality that $M$ is locally in the Monge normal form given by Proposition~\ref{prop:Euclidean_Monge_normal_form}:
$$
X\colon(\R^2,\bm{0})\to(\R^4,\bm{0}), \quad X(u,v)=(u,v,f(u,v),g(u,v)).
$$
This choice yields an adapted orthonormal frame $\{\bm{e}_1,\bm{e}_2,\bm\nu_1,\bm\nu_2\}$ of $\R^4$, where $\{\bm{e}_1,\bm{e}_2\}$ and $\{\bm\nu_1,\bm\nu_2\}$ constitute orthonormal frames of $T_{\bm{p}}M$ and $N_{\bm{p}}M$, respectively. 
Specifically, we observe that
$$
X_u(0,0)=\bm{e}_1,
\qquad
X_v(0,0)=\bm{e}_2,
$$
which implies $\text{det}[I_{\bm{p}}]=1$. 

Let the unit vectors in $S(T_{\bm{p}}M)$ be given by
$$
\bm{w}=\cos\theta\bm{e}_1+\sin\theta\bm{e}_2,\quad
\bm{w}^{\perp}=-\sin\theta\bm{e}_1+\cos\theta\bm{e}_2,\quad \theta\in\left(-\frac{\pi}2,\frac{\pi}2\right] 
$$
and those in $S(N_{\bm{p}}M)$ by
$$
\bm\nu=\cos\phi\bm\nu_1+\sin\phi\bm\nu_2,\quad
\bm\nu^{\perp}=-\sin\phi\bm\nu_1+\cos\phi\bm\nu_2,\quad \phi\in\left(-\frac{\pi}2,\frac{\pi}2\right].
$$
The projection direction $\bm{v}$ of $\pi_{\bm{v}}$ is then expressed as 
$$
\bm{v}=\cos\psi \bm{w}+\sin\psi\bm\nu^{\perp},\quad \psi\in\left(-\frac{\pi}2,\frac{\pi}2\right].
$$
By setting
$$
\bm\eta_1=-\bm{w}^{\perp},\quad
\bm\eta_2=\sin\psi \bm{w}-\cos\psi \bm\nu^{\perp},\quad
\bm\eta_3=\bm\nu,
$$
we obtain a positively oriented orthonormal frame $\{\bm\eta_1,\bm\eta_2, \bm{v}, \bm\eta_3\}$ of $\R^4$, along with the corresponding frame $\{\bm\eta_1,\bm\eta_2,\bm\eta_3\}$ of the hyperplane $\bm{v}^{\perp}$. 
Consequently, the projection $\pi_{\bm{v}}$ on $M$ along $\bm{v}$ to $\bm{v}^{\perp}$ is given by
$$
\pi_{\bm{v}}(u,v) = (\langle X(u,v),\bm\eta_1 \rangle, \langle X(u,v),\bm\eta_2 \rangle, \langle X(u,v),\bm\eta_3 \rangle).
$$

We then associate $\bm{p}\in M$ with $\pi_{\bm{v}}(0,0)\in\widetilde{M}$ and treat $\bm{p}$ as the origin of $\bm{v}^\perp$. 
In particular, $\bm{p}=\bm{0}\in \R^4$.

Suppose that $F_C$ has an $A_{\ge1}$-singularity at $(0,0)$. 
Then we set $\bm{a}=r\bm\nu$. 
Note that $\bm\nu$ constitutes a member of the frame $\{\bm\eta_1, \bm\eta_2, \bm{v}, \bm\eta_3\}$ and choose the sign of $r$ so that $\bm{a}=r\bm\nu$.


Since $\widetilde{M}$ is a regular surface at $\pi_{\bm{v}}(0,0)$, we have $\bm{v}\not\in T_{\bm{p}}M$, which is equivalent to $\sin\psi\ne0$. 
We may assume that $\sin\psi>0$ by suitably choosing the orientations of $\bm\nu^\perp$ and $\bm\eta_2$, if necessary.


The unit normal vector field to $\widetilde{M}$ is given by
$$
\widetilde{\bm{n}}=\dfrac{(\pi_{\bm{v}})_u \times (\pi_{\bm{v}})_v}{\lVert (\pi_{\bm{v}})_u \times (\pi_{\bm{v}})_v \rVert},
$$
and we thus obtain
$$
\widetilde{\bm{n}}(0,0) = (0,0,1).
$$
Therefore, we see that 
$$
\widetilde{\bm{n}}(0,0)=\bm\nu,
$$
because $\bm{v}^{\perp}$ is the subspace of $\R^4$ spanned by $\{\bm\eta_1,\bm\eta_2,\bm\eta_3\}$.

%

We consider the condition under which the distance-squared function on $\widetilde{M}$ from a point $\bm{a}\in \bm{v}^{\perp}\subset \R^4$ given by
$$
d_{\bm{a}}\colon (\R^2,\bm{0})\to\R,\quad d_{\bm{a}}(u,v) = \lVert \pi_{\bm{v}}(u,v)-\bm{a}\rVert^2
$$
has an $A_{\geq 1}$-singularity. 

It follows that $d_{\bm{a}}$ has an $A_{\ge1}$-singularity at $(0,0)$ if and only if
$$
\langle (\pi_{\bm{v}})_u(0,0), -\bm{a}\rangle = \langle (\pi_{\bm{v}})_v(0,0), -\bm{a}\rangle = 0
$$
which is equivalent to 
\begin{equation*}
\begin{array}{ccccc}
 & \bm{v}^\perp & \subset & \R^4 & \\ [-1mm]
 & \rotatebox[origin=c]{90}{$\in$}
 & & \rotatebox[origin=c]{90}{$\in$} & \\ [-1mm]
\bm{a}= & r\widetilde{\bm{n}}(0,0)
& = & r\bm\nu,
\end{array}
\end{equation*}
where $|r|$ is the radius of the sphere 
$$
P_{\bm{v}}(C_{\bm{v},\bm{a},r})=\{P_{\bm{v}}(\bm{x})\in \bm{v}^\perp \colon\lVert P_{\bm{v}}(\bm{x})-\bm{a}\rVert^2=r^2, \bm{x}\in \R^4\}.
$$ 
Hence, $F_C$ has an $A_{\geq 1}$-singularity at $(0,0)$ if and only if $d_{\bm{a}}$ has an $A_{\geq 1}$-singularity there.

Since
$$
F_C(u,v)=\lVert \pi_{\bm{v}}(u,v)-\bm{a}\rVert^2=d_{\bm{a}}(u,v),
$$
Lemma~\ref{lem:main} now yields this theorem.

\end{proof}

Theorem~\ref{thm:main} (3) implies that the axis of $A_{\geq 3}$-contact spherical cylinder $C_{\bm{v},\bm{a},r}$ with a regular surface in $\R^4$ is characterized as a higher-dimensional extension of the axis of $A_{\geq 3}$-contact cylinder of a regular surface in $\R^3$. 
Therefore, the notion of cylindrical directions of regular surfaces in $\R^3$ can be extended to regular surfaces in $\R^4$ as follows: 
Assume that a projected surface $\widetilde{M}=P_{\bm{v}}(M)$ is regular at $P_{\bm{v}}(\bm{p})$. 
Then the direction $\bm{v}$ of the axis of $C_{\bm{v},\bm{a},r}$ is called a \textit{cylindrical direction of a regular surface $M$ in $\R^4$ at $\bm{p}$ derived from its contact with a spherical cylinder $C_{\bm{v},\bm{a},r}$} if $M$ has an $A_{\geq 3}$-contact with $C_{\bm{v},\bm{a},r}$ at $\bm{p}$. 
We simply call it a spherical-cylindrical direction of $M$. 

\begin{remark}
\begin{enumerate}
\item
The definition of a cylindrical direction of a regular surface $M_3$ in $\R^3$ requires that the apparent contour of $M_3$ along the direction of the axis of a cylinder having degenerate contact with the surface have no inflection point at the point corresponding to $\bm{p}\in M_3$, since the apparent contour has no finite focal point at an inflection point.

On the other hand, the definition of a spherical-cylindrical direction of a surface $M$ in $\R^4$ does not require that the projected regular surface of $M$ along the direction of the axis of a spherical cylinder having degenerate contact with the surface be non-parabolic at the point corresponding to $\bm{p}\in M$. Indeed, even at a parabolic point, although one principal curvature vanishes, the projected regular surface has a finite focal point associated with the nonzero principal curvature.

\item
There exist at most three cylindrical directions at a non-parabolic point of a regular surface in $\R^3$ (see \cite{FHN2017}). 
In more detail, there exist generically one or three cylindrical direction at a non-parabolic point of the surface. 
In \cite{TSH2025}, $\Delta_3$-curve, which is the locus of points where there exist two cylindrical directions of the surface, was introduced. 
Moreover, they introduced the $V(2)^{ac}$-curve which is the locus of points where a cylinder has an $A_{\geq 4}$-contact with the surface. 
The $V(2)^{ac}$-curve is also the locus of points where the apparent contour of the surface associated with the orthogonal projection along the cylindrical direction has a degenerate vertex. 

However,  in this paper, we do not deal with the number of the spherical-cylindrical directions of a regular surface in $\R^4$. 
Similarly, a higher-dimensional extension of $\Delta_3$-curve and those of $V(2)^{ac}$-curve are not addressed in this paper. 
These remain the subject of future research. 

\item
In \cite{DSS2026}, inflections and vertices of the apparent contour of a generic regular surface in $\R^4$ associated with the parallel projection $\bm{P}_{\pi}$ along the 2-planes $\pi$ onto transverse planes were studied. 
The contact between the surface and 2-planes $\pi$ is characterized by the singularities of $\bm{P}_{\pi}$. 
Therefore, we expect a relationship between the vertices of the apparent contour of the surface in $\R^4$ and their contact with toroidal cylinders $S^1(r) \times \R^2$.
\end{enumerate}
\end{remark}

Under the same setting as in the proof of Theorem~\ref{thm:main}, suppose that $F_C$ has an $A_{\geq 1}$-singularity at $(0,0)$. 
Then the Hessian matrix of $F_C$ at $(0,0)$ is given by
\begin{equation*}
H(F_C)(0,0)
=2\begin{pmatrix}
c_{20} & c_{11}\\
c_{11} & c_{02}
\end{pmatrix},
\end{equation*}
where
\begin{align*}
c_{20} &=\sin^2\theta+\cos^2\theta\sin^2\psi -r\bigl((h_1+\lambda)\cos\phi+h_2\sin\phi\bigr),\\
c_{11} &=-\cos\theta\sin\theta\cos^2\psi-r\mu\sin\phi,\\
c_{02} &=\cos^2\theta+\sin^2\theta\sin^2\psi -r\bigl((h_1-\lambda)\cos\phi+h_2\sin\phi\bigr).
\end{align*}
Moreover,
\begin{equation*}
\det H(F_C)(0,0) =4\bigl(A(\psi)+rB(\theta,\phi,\psi)+r^2Q(\phi)\bigr),
\end{equation*}
where
\begin{align*}
A(\psi) &=\sin^2\psi,\quad 
B(\theta,\phi,\psi) =-\cos^2\psi\langle II_{\bm{p}}(\bm{w},\bm{w}),\bm\nu\rangle -2\sin^2\psi\langle \bm{H},\bm\nu\rangle,\\
Q(\phi) &=(h_1^2-\lambda^2)\cos^2\phi +2h_1h_2\cos\phi\sin\phi +(h_2^2-\mu^2)\sin^2\phi=\det[II_{\bm{p},\bm\nu}].
\end{align*}

We now consider the geometrical meaning of $Q(\phi)$. 
The function $Q(\phi)$ is given by 
$$
Q(\phi) = \text{det}\left[\widetilde{II}_{\bm{p}}\right],
$$
where $\widetilde{II}_{\bm{p}}$ is the second fundamental form of a projected regular surface $P_{\bm{v}}(M)=\widetilde{M}$ (see \cite{Hasegawa2026}). 
It follows that $Q(\phi)=0$ if and only if $P_{\bm{v}}(\bm{p})$ is either a parabolic point or a flat umbilic point of $\widetilde{M}$.
On the other hand, from a viewpoint of $M$, it follows from $Q(\phi)=\det[II_{\bm{p},\bm\nu}]$ that $Q(\phi)=0$ if and only if there exist asymptotic directions $\bm{w} \in T_pM$ and binormal directions $\bm\nu \in N_pM$ associated with $\bm{w}$ at $\bm{p}$.
Moreover, let $\bm{p}$ be a non-inflection point. 
Since $\text{disc}(\text{det}[II_{\bm{p},\bm\nu}])=\delta(\bm{p})$, $Q(\phi)=0$ has no real roots, one repeated real root or distinct two roots if and only if $\delta(\bm{p})$ is negative, zero, positive which are equivalent to $\bm{p}$ being an elliptic point, parabolic point or hyperbolic point, respectively. 



\begin{thebibliography}{99}

\bibitem{AOS2026}
F.~Alharbi and R.~Oset Sinha, 
Height functions on a cuspidal edge with a distinguished singular curve, \href{https://doi.org/10.4171/pm/2149}
{Port. Math. \textbf{83} (2026), 167--185.}

\bibitem{Arnold1986}
V.~I.~Arnol'd, S.~M.~Guse{\u\i}n-Zade and A.~N.~Varchenko,
Singularities of differentiable maps.\ Vol.~I.
\href{https://doi.org/10.1007/978-0-8176-8340-5}
{Monogr.\ Math. 82, (Birkh$\ddot{\text{a}}$user Boston, Inc., 1986.)}





\bibitem{BROS2019}
P.~Benedini Riul and R.~Oset Sinha, 
A relation between the curvature ellipse and the curvature parabola, 
\href{https://doi.org/10.1515/advgeom-2019-0002}
{Adv. Geom. \textbf{19} (2019) 389--399.}

\bibitem{BGT1995}
J.~W.~Bruce, P.~J.~Giblin and F.~Tari,
Families of surfaces: height functions, Gauss maps and duals, 
Pitman Res. Notes Math. \textbf{333} (1995), 148--178. 

\bibitem{BGT1999}
J.~W. Bruce, P.~J. Giblin and F.~Tari, 
Families of surfaces: Focal sets, ridges and umbilics, 
\href{https://doi.org/10.1017/S0305004198003004}
{Math. Proc. Cambridge Philos. Soc. \textbf{125} (1999), 243--268.} 

\bibitem{BN1998}
J.~W.~Bruce and A.~C.~Nogueira,
Surfaces in $\R^4$ and duality, 
\href{https://doi.org/10.1093/qmathj/49.4.433}
{Quart. J. Math. Oxford \textbf{49} (1998), 433-443. }

\bibitem{BT2002}
J.~W.~Bruce and F.~Tari,
Families of surfaces in $\R^4$,
\href{https://doi.org/10.1017/S0013091500000213}
{Proc. Edinb. Math. Soc. \textbf{45} (2002), 181--203.}

\bibitem{David1983}
J.~M.~S.~David, 
Projection-generic curves,
\href{https://doi.org/10.1112/jlms/s2-27.3.552}
{J. London Math. Soc. (2) \textbf{27} (1983), 552--562. }

\bibitem{DSilva2021}
J.~L.~Deolindo-Silva, 
Cross-ratio invariants for surfaces in 4-Space, 
\href{https://doi.org/10.1007/s00574-020-00221-w}
{Bull. Braz. Math. Soc. New Series \textbf{52} (2021), 591--612.} 


\bibitem{DSK2019}
J.~L.~Deolindo-Silva and Y.~Kabata, 
Projective classification of jets of surfaces in 4-space, 
\href{https://doi.org/10.32917/hmj/1554516037}
{Hiroshima Math. J. \textbf{49} (2019), 35--46.}

\bibitem{DSS2026}
J.~L.~Deolindo-Silva and M.~Salarinoghabi, 
On the contact of surfaces in 4-space with 2-planes and their apparent contours, 
\href{https://doi.org/10.1007/s13163-026-00603-6}
{Rev. Mat. Complut. (2026). }


\bibitem{Francisco2023}
A.~P.~Francisco, 
Functions on a swallowtail, 
\href{https://doi.org/10.21711/231766362023/rmc534}
{Mat. Contemp. \textbf{53} (2023), 5--74. }

\bibitem{FH2012-1}
T.~Fukui and M.~Hasegawa,
Singularities of parallel surfaces, 
\href{https://doi.org/10.2748/tmj/1347369369}
{Tohoku Math. J. \textbf{64} (2012), 387--408.}

\bibitem{FH2012-2}
T.~Fukui and M.~Hasegawa,
Fronts of Whitney umbrella-- a differential geometric approach via blowing up,
\href{http://dx.doi.org/10.5427/jsing.2012.4c}
{J. Singul., {\bf 4} (2012), 35--67. }

\bibitem{FH2013}
T.~Fukui and M.~Hasegawa,
Height functions on Whitney umbrellas,
\href{http://hdl.handle.net/2433/207802}
{K\^oky\^uroku Bessatsu, {\bf B38} (2013), 153--168.}

\bibitem{FHN2017}
T.~Fukui, M.~Hasegawa and K.~Nakagawa, 
Contact of a regular surface in Euclidean $3$-space with cylinders and cubic binary differential equations,
\href{https://doi.org/10.2969/jmsj/06920819}
{J. Math. Soc. Japan \textbf{69} (2017), 819--847.}

\bibitem{FH2025-1}
T.~Fukui and M.~Hasegawa, 
Height functions on singular surfaces parameterized by smooth maps $\mathcal{A}$-equivalent to $S_k$, $B_k$, $C_k$ and $F_4$, 
\href{https://doi.org/10.1007/s40687-025-00565-z}
{Res. Math. Sci., \textbf{12} (2025). }

\bibitem{FH2025-2}
T.~Fukui and M.~Hasegawa, 
Distance-squared functions on singular surfaces parameterized by smooth maps $\mathcal{A}$-equivalent to  $S_k$, $B_k$, $C_k$ and $F_4$, 
\href{https://doi.org/10.2969/aspm/08910113}
{Adv. Stud. Pure Math., \textbf{89} (2025), 113--151.}

\bibitem{GMRFR2000}
R. Garcia, D. K. H. Mochida, M. C. Romero-Fuster and M. A. S. Ruas,
Inflection points and topology of surfaces in 4-space, 
\href{https://doi.org/10.2307/118196}
{Trans. Am. Math. Soc. \textbf{352} (2000), 3029--3043.}

\bibitem{Hasegawa2026}
M.~Hasegawa, 
A curvature-ellipse-based Monge normal form of a regular surface in $\R^4$ and its applications, 
\href{https://doi.org/10.48550/arXiv.2610.10132}{preprint.} 

\bibitem{HS2025}
M.~Hasegawa and K.~Saji, 
On cylindrical directions of cross caps and cuspidal edges, 
\href{https://doi.org/10.14492/hokmj/2024-810}
{Hokkaido Math. J. \textbf{54} (2025), 425--446.}




\bibitem{IFRT2015}
S.~Izumiya, M.~C.~Romero Fuster, M.~A.~S.~Ruas and F.~Tari, Differential Geometry from a Singularity Theory Viewpoint,
\href{https://doi.org/10.1142/9108}
{World Scientific, 2015.}

\bibitem{Little1969}
J.~A. Little,
On singularities of submanifolds of higher dimensional Euclidean spaces, 
\href{https://doi.org/10.1007/BF02411172}
{Ann. Mat. Pura Appl. \textbf{83} (1969), 261--335.}


\bibitem{Mather1968}
J.~N. Mather, 
Stability of $C^\infty$ mappings. III. Finitely determined map-germs, 
\href{https://doi.org/10.1007/BF02684889}
{Publ. Math. Inst. Hautes \'Etudes Sci. \textbf{35} (1968), 279--310. }

\bibitem{MRFR1995}
D.~K.~H.~Mochida, M.~D.~C.~Romero Fuster and M.~A.~S.~Ruas,
The geometry of surfaces in 4-space from a contact viewpoint, 
\href{https://doi.org/10.1007/BF01265348}
{Geom. Dedicata \textbf{54} (1995), 323--332.}

\bibitem{Montaldi1986-1}
J.~A.~Montaldi,
On contact between submanifolds, 
\href{https://doi.org/10.1307/mmj/1029003348}
{Michigan Math. J. {\bf 33} (1986), 81--85.}

\bibitem{Mond1982}
D.~M.~Q.~Mond, 
Classification of certain singularities and applications to differential geometry,
PhD thesis, University of Liverpool (1982).

\bibitem{Mond1985}
D.~M.~Q.~Mond, 
On the classification of germs of maps from $\R^2$ to $\R^3$, 
\href{https://doi.org/10.1112/plms/s3-50.2.333}
{Proc. Lond. Math. Soc. \textbf{50} (1985), 333--369.}


\bibitem{NBT2007}
J.~J.~Nu\~no-Ballesteros and F.~Tari, 
Surfaces in $\R^4$ and their projections to 3-spaces, 
\href{https://doi.org/10.1017/S0308210506000758}
{Proc. R. Soc. Edinburgh, Sect. A \textbf{137} (2007), 1313--1328.}

\bibitem{OSS2018}
R. Oset Sinha and K. Saji, 
On the geometry of folded cuspidal edges, 
\href{https://doi.org/10.1007/s13163-018-0257-6}
{Rev. Mat. Complut. \textbf{31} (2018), 62--650. }

\bibitem{OST2015}
R.~Oset Sinha and F.~Tari, 
Projections of surfaces in $\R^4$ to $\R^3$ and the geometry of their singular images, 
\href{https://doi.org/10.4171/RMI/825}
{Rev. Mat. Iberoam. \textbf{31} (2015), 33--50.}

\bibitem{OST2018}
R.~Oset Sinha and F.~Tari, 
On the flat geometry of the cuspidal edge, 
\href{https://doi.org/10.18910/70079}
{Osaka J. Math. \textbf{55} (2018), 393--421. }

\bibitem{Porteous1971}
I.~R.~Porteous, 
The normal singularities of a submanifold, 
\href{https://doi.org/10.4310/jdg/1214430015}
{J. Differ. Geom. \textbf{5} (1971), 543--564. }

\bibitem{Porteous1983}
I.~R.~Porteous, 
The normal singularities of surfaces in $\R^3$, 
Singularities, Part 2 (Arcata, Calif., 1981), Proc. Sympos. Pure Math., \textbf{40}, Amer. Math. Soc., Providence, RI, 1983, 379--393.

\bibitem{RFSC2013}
M.~C.~Romero-Fuster and E.~Sanabria-Codesal, 
Conformal invariants and spherical contacts of surfaces in $\R^4$, 
\href{https://doi.org/10.1007/s13163-011-0086-3}
{Rev. Mat. Complut. \textbf{26} (2013), 215--240. }

\bibitem{Segre1938}
B.~Segre, 
Famiglie di ipersuperficie isoparametriche negli spazi euclidei ad un qualunque numero di dimensioni,
Atti Accad. Naz. dei Lincei Cl. Sci. Fis. Mat. Nat. Rend. Lincei Mat. Appl. \textbf{27} (1938), 203--207.

\bibitem{Takahashi1970}
T.~Takahashi, 
Homogeneous hypersurfaces in space of constant curvature, 
\href{https://doi.org/10.2969/jmsj/02230395}
{J. Math. Soc. Japan, \textbf{22} (1970), 395--410.}


\bibitem{TSH2025}
F.~Tari, M.~Salarinoghabi and M.~Hasegawa, 
Geometric deformations of discriminants and apparent contours, 
\href{https://doi.org/10.1007/978-3-031-87016-3}
{Springer, New York (2025). }

\bibitem{Teramoto2026}
K.~Teramoto, 
Height Functions on Generalized Cuspidal Edges in Certain Directions, 
\href{https://doi.org/10.1007/s00025-026-02305-w}
{Results Math. \textbf{79} (2026).}

\bibitem{Thom1972}
R.~Thom, 
Stabilit{\'e} structurelle et morphog{\'e}n{\`e}se: Essai d'une th{\'e}orie g{\'e}n{\'e}rale des mod{\`e}les, 
\href{https://archive.org/details/stabilitestructu0000thom_y5k9}
{W. A. Benjamin, Reading, Mass., 1972.}


\end{thebibliography}
\end{document}